\documentclass[11pt]{article}
\usepackage[a4paper, tmargin=3cm, bmargin=3.0cm, lmargin=2.5cm, rmargin=2.5cm, textheight=24cm, textwidth=16cm]{geometry}
\usepackage{setspace}
\usepackage{amsfonts}
\usepackage{epsfig}
\usepackage{latexsym}
\usepackage{amsthm}
\usepackage{amsmath}
\usepackage{amssymb}
\usepackage{array, xcolor}
\usepackage{enumerate}
\usepackage{enumitem}
\usepackage{graphicx}
\usepackage{breqn}
\usepackage{thmtools, thm-restate}

\newtheorem{theorem}{Theorem}[section]

\newtheorem{example}{Example}[section]

\newtheorem{corollary}{Corollary}[section]

\newtheorem{lemma}{Lemma}[section]

\newlist{notes}{enumerate}{1}
\setlist[notes]{label=Note:,leftmargin=*}

\title{Spectral radius and second largest eigenvalues of power graphs of finite groups}
\author{ Priti Prasanna Mondal$^{1, 2, *}$, Basit Auyoob Mir$^{2}$, Fouzul Atik$^{2}$
 \\
 \small $^{1}$ Department of Mathematics, Indian Institute of Technology Hyderabad, Hyderabad 502284, India\\
 \small$^{2}$ Department of Mathematics, SRM University-AP, Andhra Pradesh 522240, India.\\ \small e-mail: pritiprasanna1992@gmail.com, mirbasit553@gmail.com, fouzulatik@gmail.com\\
 \small $^{*}$Corresponding author: \tt{pritiprasanna1992@gmail.com} \\
 }

\date{}

\begin{document}
\maketitle{}
\begin{abstract}
Consider a group $\mathbb{G}$ and construct its power graph, whose vertex set consists of the elements of $\mathbb{G}$. Two distinct vertices (elements) are adjacent in the graph if and only if one element can be expressed as an integral power of the other. In this article, we improved the bounds of the spectral radius of the power graphs of the cyclic group $C_{n}$, the dihedral group $\mathcal{D}_{2n}$, and the dicyclic group $\mathcal{Q}_{4n}$. For $n\neq p^{m},$ the power graph of the cyclic group $C_{n}$ is not a complete multipartite graph. We find the second largest eigenvalue bounds of the same with the clique number. In some cases, we find the bounds are exact if and only if they belong to a particular family of graphs. Lastly, we work on the distance spectral radius of the power graphs of the same groups\\

\end{abstract}

\section{Introduction}
Algebraic graph theory is a trending research area in the field of mathematics. Recently, much work has been done in this area to study algebraic structures through the graphical representation of semigroups, groups, rings, or vector spaces. In this area, the power graph of groups is a research topic. The power graph $\mathcal{P}(\mathbb{G})$ of a group $\mathbb{G}$ is an undirected graph where every element of the group $\mathbb{G}$ represents a vertex, and two vertices $a, b \in \mathbb{G}$ are adjacent if and only if $a \neq b$ and $a^m = b$ or $b^m = a$ for some $m \in \mathbb{N}$. The concept of directed power graphs was first introduced by Kelarev and Quinn \cite{Kela}. Motivated by this, Chakrabarty et al. \cite{Chak} introduced the concept of an undirected power graph $\mathcal{P}(\mathbb{G})$ of a group $\mathbb{G}$ and proved that the power graph $\mathcal{P}(\mathbb{G})$ is complete if and only if $\mathbb{G}$ is a cyclic group of order $p^m$ for some prime $p$ and non-negative integer $m$. Subsequently, researchers have studied the power graph from a graph-theoretic perspective in \cite{Came, srip2, srip3, Doos}. Later on, the spectral properties of the power graph have also been examined.

\subsection{Groups and power graphs}

Let $\mathbb{G}$ be a finite group. The \emph{power graph} $\mathcal{P}(\mathbb{G})$ is the undirected graph whose vertex set is $\mathbb{G}$, where two distinct vertices $x,y\in\mathbb{G}$ are adjacent if and only if one is an integral power of the other.

In this paper, we mainly consider the following families of groups:
\begin{itemize}
\item the cyclic group $C_n=\langle g\rangle$ of order $n$;
\item the dihedral group $\mathcal{D}_{2n}=\langle a,b \mid a^n=e,\; b^2=e,\; ba=a^{-1}b\rangle$ of order $2n$;
\item the dicyclic group $\mathcal{Q}_{4n}=\langle a,b \mid a^{2n}=e,\; a^n=b^2,\; ab=ba^{-1}\rangle$ of order $4n$.
\end{itemize}

For a graph $G$, the clique number is denoted by $\omega(G)$ and represents the order of a largest complete subgraph of $G$.

\subsection{Graphs and spectra}

Let $G=(V(G),E(G))$ be a simple graph with vertex set
$V(G)=\{v_1,v_2,\dots,v_n\}$. The degree of a vertex $v_i$ is denoted by $d_i$, and the average degree of $G$ is defined as
\[
d_{\mathrm{avg}}=\frac{1}{n}\sum_{i=1}^n d_i.
\]

The adjacency matrix of $G$ is the $n\times n$ matrix
$A(G)=(a_{ij})$, where $a_{ij}=1$ if $v_i$ is adjacent to $v_j$ and
$a_{ij}=0$ otherwise. The matrix $A(G)$ is real and symmetric, and its eigenvalues are arranged as
\[
\lambda_1(G)\ge \lambda_2(G)\ge \cdots \ge \lambda_n(G).
\]
The largest eigenvalue $\lambda_1(G)$ is called the \emph{spectral radius} of $G$, while $\lambda_2(G)$ is referred to as the \emph{second largest eigenvalue}.

The eigenvalues of graphs have wide applications in the fields of science and engineering. The spectral radius is crucial for modeling virus propagation in computer networks. Moreover, the eigenvalues of the adjacency matrix serve as valuable tools for protecting personal data privacy in certain databases. The spectral radius is linked to how quickly information or viruses propagate through a network. It is also used to identify influential nodes within a network. Computing the spectrum of graphs is an interesting research field; see \cite{Baba, srip3, Brow, srip}.

\subsection{Distance matrix and distance spectrum}

Let $G$ be a connected graph. The distance matrix of $G$ is the $n\times n$ matrix $D(G)=(d_{ij})$, where $d_{ij}$ denotes the length of a shortest path between $v_i$ and $v_j$. The transmission of a vertex $v_i$ is defined by
\[
\operatorname{Tr}(v_i)=\sum_{j=1}^n d_{ij}.
\]
We denote by $\operatorname{Tr}_{\min}$, $\operatorname{Tr}_{\max}$, and
$\operatorname{Tr}_{\mathrm{avg}}$ the minimum, maximum, and average transmission of $G$, respectively.

The distance matrix $D(G)$ is real and symmetric, and its eigenvalues are ordered as
\[
\rho_1(G) > \rho_2(G) \ge \cdots \ge \rho_n(G),
\]
where $\rho_1(G)$ is called the \emph{distance spectral radius} of $G$.

The distance spectra of graphs are used extensively in the fields of science and engineering. The distance spectral radius is a useful molecular descriptor in QSPR (Quantitative Structure-Property Relationship) modeling. Zhou and Trinajstic provided upper and lower bounds for the distance spectral radius in terms of the number of vertices, the Wiener index, and the Zagreb index in \cite{zhou}. In \cite{Das}, Das determined the upper and lower bounds for the distance spectral radius of connected bipartite graphs and characterized the graphs for which these bounds are exact. In \cite{Ind}, Indulal found sharp bounds on the distance spectral radius and the distance energy of graphs. Atik et al. \cite{Atik2} identified all the distance eigenvalues of $m$-generation $n$-prism graphs. In \cite{Atik}, the authors gave distance eigenvalue bounds other than the spectral radius using the smallest Ger$\check{\text{s}}$gorin disc. Pirzada et al. obtained the distance signless Laplacian spectrum of some algebraic graphs in \cite{pirzada}.


\subsection{Linear algebraic results}

 We recall several well-known results concerning quotient matrices, equitable partitions, and the spectral radius of non-negative matrices, which will be used in the sequel.

Let $A$ be an $n\times n$ matrix whose rows and columns are indexed by the set
$X=\{1,2,\ldots,n\}$. Consider a partition
\[
\pi=\{X_1,X_2,\ldots,X_m\}
\]
of the index set $X$. The \emph{characteristic matrix} of the partition $\pi$ is the $n\times m$ matrix
$C=(c_{ij})$, where
\[
c_{ij}=
\begin{cases}
1, & \text{if } i\in X_j,\\
0, & \text{otherwise}.
\end{cases}
\]

According to the partition $\pi$, the matrix $A$ can be written in block form as
\[
A=
\begin{pmatrix}
A_{11} & A_{12} & \cdots & A_{1m}\\
A_{21} & A_{22} & \cdots & A_{2m}\\
\vdots & \vdots & \ddots & \vdots\\
A_{m1} & A_{m2} & \cdots & A_{mm}
\end{pmatrix},
\]
where $A_{ij}=A[X_i:X_j]$ for $1\le i,j\le m$.
If $q_{ij}$ denotes the average row sum of the block $A_{ij}$, then the matrix
$Q=(q_{ij})$ is called the \emph{quotient matrix} of $A$ corresponding to the partition $\pi$.
If each block $A_{ij}$ has constant row sum, then the partition $\pi$ is said to be
\emph{equitable}.

The following classical result describes the relation between the eigenvalues of a matrix and those of its quotient matrix.

\begin{theorem}\label{quotient}
\emph{\cite{Brow}}
Let $A$ be a symmetric matrix whose rows and columns are partitioned according to
$\{X_1,X_2,\ldots,X_m\}$, and let $Q$ be the corresponding quotient matrix.
Then the eigenvalues of $Q$ interlace the eigenvalues of $A$.
\end{theorem}

For equitable partitions, the spectral radius of a matrix is preserved.

\begin{theorem}\label{spec}
\emph{\cite{Atik3}}
Let $A$ be a non-negative matrix and let $Q$ be a quotient matrix corresponding to an equitable partition of $A$. Then $A$ and $Q$ have the same spectral radius.
\end{theorem}

We also recall the Perron-Frobenius theorem, which provides fundamental properties of the spectral radius of non-negative matrices.

\begin{theorem}[Perron--Frobenius]\label{perron}
\emph{\cite{horn}}
Let $A\in M_n$ be an irreducible non-negative matrix with $n\ge 2$. Then:
\begin{enumerate}
\item $\rho(A)>0$;
\item $\rho(A)$ is a simple eigenvalue of $A$;
\item there exists a unique positive eigenvector $x=(x_i)$ corresponding to $\rho(A)$,  such that
\[
\sum_{i=1}^n x_i = 1.
\]
\end{enumerate}
\end{theorem}

The paper is organized as follows. In Section~1, we collect necessary definitions and preliminary results. Section~2 is devoted to bounds on the spectral radius and the second largest eigenvalue of power graphs, respectively. In Section~3, we investigate the distance spectral radius of power graphs of cyclic, dihedral, and dicyclic groups.

Throughout the paper, we follow standard terminology and notation from algebraic graph theory. Additional notation will be introduced when needed.

\section{Main Results}
\subsection{Cyclic Group}
In \cite{srip}, Chattopadhyay et al. provided bounds on the spectral radius of the power graph of the cyclic group using the maximum degree and the minimum degree of the graph. Here, we give an improved lower bound on the spectral radius by using the average degree of the graph and proof by the quotient matrix technique. 
\begin{theorem}\label{spec_cn}
Let $\lambda_{1}(\mathcal{P}(C_{n}))$ denote the spectral radius of $\mathcal{P}(C_{n})$. For $n\geq 3$, the improved spectral radius (lower) bound is $$ \lambda_{1}(\mathcal{P}(C_{n})) \ge \frac{1}{2}\left[(d_{avg}-1) + \sqrt{(d_{avg}+1-2l)^{2}+4l(n-l)}\right].$$
Where $l=\phi(n)+1$ and $d_{avg}$ is the average degree of the non-identity non-generator elements of $C_{n}$. Moreover, equality holds if $\mathcal{P}(C_{n})$ is a complete graph. 
\end{theorem}
\begin{proof}
    Let $V_{0},~V_{1}$ and $V_{2}$ be three partition sets of the vertices of the power graph of the cyclic group $C_{n}.$ The set $V_{0}$ contains the identity element of $C_{n}$, the set $V_{1}$ contains all generators of $C_{n}$, and the set $V_{2}$ contains the rest of the vertices. Then $|V_{1}|= 1,~ |V_{1}|= \phi(n)=l-1$ (say), where $\phi(n)$ is Euler's $\phi$ function. We ordered the sets $V_0,~V_1$ and $V_2$ as $\{v_{0}=e\},~\{ v_{2}, \cdots, v_{l}\}$ and $\{ v_{l+1}, v_{l+2}, \cdots, v_{n}\}$ respectively. Then $A(\mathcal{P}(C_{n}))$ with rows and columns indexed by $\{ V_{0}, V_{1}, V_{2}\}$ can be written as 
     \begin{eqnarray}\label{APC}\displaystyle {A(\mathcal{P}(C_{n}))}&=&{\begin{bmatrix}  {0}&{{\textbf{1}}^{T}_{l-1}} &{\textbf{1}^{T}_{n-l}}\\
  {\textbf{1}_{l-1}} &  {J_{l-1}-I_{l-1}} &{J_{(l-1)\times (n-l)}}\\
 {\textbf{1}_{n-l}}&{J_{(n-l)\times l}} &{A(\mathcal{P}(V_{2}))_{(n-l)}} 
 \end{bmatrix}}.
 \end{eqnarray}
Where $J_{m \times n}$ is a matrix with all entries $1$, $\textbf{1}_{k}$ is a vector of size $k$ with all entries $1$,  $I$ is the identity matrix, and $A(\mathcal{P}(V_{2}))=(a_{ij})$ is the adjacency matrix of the power graph induced by the vertices of $V_{2}$. Now the quotient matrix of the $A(\mathcal{P}(C_{n}))$ matrix with respect to the partition $\pi=\{ V_{0}, V_{1}, V_{2}\}$ is 
 \begin{eqnarray*}\displaystyle {Q_{C_{n}}}&=&{\begin{bmatrix}  {0}&{l-1} &{n-l}\\
 {1}&{l-2} &{n-l}\\
 {1}&{l-1} &{d_{avg}-l} 
 \end{bmatrix}}.
 \end{eqnarray*}
 The spectral radius of $Q_{C_{n}}$ is $\mu_{1}({Q_{C_{n}}})=\frac{1}{2}\left[(d_{avg}-1) + \sqrt{(d_{avg}+1-2l)^{2}+4l(n-l)}\right].$
So, by the Hammer interlacing of eigenvalues ( see the Theorem \ref{quotient}), we get $\lambda_{1}(\mathcal{P}(C_{n}))\ge \mu_{1}({Q_{C_{n}}}).$ \\ Therefore, $\lambda_{1}(\mathcal{P}(C_{n})) \ge \frac{1}{2}\left[(d_{avg}-1) + \sqrt{(d_{avg}+1-2l)^{2}+4l(n-l)}\right]$.

If $\mathcal{P}(C_{n})$ is a complete graph, then $d_{avg}=n-1$. For a complete graph $\lambda_{1}(\mathcal{P}(C_{n}))=n-1$ and it is easy to verify that $\frac{1}{2}\left[(d_{avg}-1) + \sqrt{(d_{avg}+1-2l)^{2}+4l(n-l)}\right]=n-1$. So, equality holds if it $\mathcal{P}(C_{n})$ is a complete graph.
\end{proof}

\textbf{Comparison:} It is clear that $d_{avg}\ge d_{min}$; the equality holds for regular graphs only.
In \cite{srip}, the given spectral radius lower bound is $\lambda_{1}(\mathcal{P}(C_{n})) \ge \frac{1}{2}\left[(d_{min}-1) + \sqrt{(d_{min}+1-2l)^{2}+4l(n-l)}\right]$. And our result is $\lambda_{1}(\mathcal{P}(C_{n})) \ge \frac{1}{2}\left[(d_{avg}-1) + \sqrt{(d_{avg}+1-2l)^{2}+4l(n-l)}\right]$. These two lower bounds are comparable and our bound is better, as we have $$\frac{1}{2}\left[(d_{avg}-1) + \sqrt{(d_{avg}+1-2l)^{2}+4l(n-l)}\right]\ge \frac{1}{2}\left[(d_{min}-1) + \sqrt{(d_{min}+1-2l)^{2}+4l(n-l)}\right].$$

Now one immediate result on the second largest eigenvalue of $\mathcal{P}(C_{n})$ can be stated as follows:
\begin{theorem}\label{2nd_spec_cn}
Let $\lambda_{2}(\mathcal{P}(C_{n}))$ be the second largest eigenvalue of $\mathcal{P}(C_{n})$. For $n\geq 3$, the second largest eigenvalue is $\lambda_{2}(\mathcal{P}(C_{n})) \ge -1.$
The equality holds if and only if $\mathcal{P}(C_{n})$ is a complete graph. 
\end{theorem}
\begin{proof}
The quotient matrix of the $A(\mathcal{P}(C_{n}))$ matrix with respect to the partition $\pi=\{ V_{0}, V_{1}, V_{2}\}$ (same as in Theorem {\ref{spec_cn}}) is 
 \begin{eqnarray*}\displaystyle {Q_{C_{n}}}&=&{\begin{bmatrix}  {0}&{l-1} &{n-l}\\
 {1}&{l-2} &{n-l}\\
 {1}&{l-1} &{d_{avg}-l}
 \end{bmatrix}}.
 \end{eqnarray*}
 The eigenvalues of $Q_{C_{n}}$ is $\mu_{1}({Q_{C_{n}}})=\frac{1}{2}\left[(d_{avg}-1) + \sqrt{(d_{avg}+1-2l)^{2}+4l(n-l)}\right],~ \mu_{2}({Q_{C_{n}}})=-1$ and $\mu_{3}({Q_{C_{n}}})=\frac{1}{2}\left[(d_{avg}-1) - \sqrt{(d_{avg}+1-2l)^{2}+4l(n-l)}\right].$ It is clear that $\mu_{1}({Q_{C_{n}}})\ge\mu_{2}({Q_{C_{n}}})\ge\mu_{3}({Q_{C_{n}}}).$
So, by the Hammer interlacing of eigenvalues (see Lemma \ref{quotient}), we get $\lambda_{2}(\mathcal{P}(C_{n}))\ge \mu_{2}({Q_{C_{n}}}).$  Therefore, $\lambda_{2}(\mathcal{P}(C_{n})) \ge -1.$ 

If $\mathcal{P}(C_{n})$ is a complete graph, then $\lambda_{2}(\mathcal{P}(C_{n}))=-1.$ So, the equality holds if $\mathcal{P}(C_{n})$ is a complete graph.

\end{proof}

We can say something more about second largest eigenvalue of $\mathcal{P}(C_{n}),$ if we remove the class of cyclic group of order $n=p^m.$

\begin{theorem}\label{multipart}
    Let $\lambda_{2}(\mathcal{P}(C_{n}))$ be the second largest eigenvalue of $\mathcal{P}(C_{n})$. For $n\neq p^m$, the second largest eigenvalue is $\lambda_{2}(\mathcal{P}(C_{n})) > 0.$
\end{theorem}
To prove this, it is sufficient to show that the power graph $\mathcal{P}(C_{n})$ is not complete multipartite for $n\neq p^m.$ We prove this in the following lemma.
\begin{lemma}
    The power graph $\mathcal{P}(C_{n})$ of cyclic group $C_{n}$ is not complete multipartite for $n\neq p^m.$
\end{lemma}
\begin{proof}
Let $n=p_{1}^{\alpha_{1}}p_{2}^{\alpha_{2}}\cdots p_{r}^{\alpha_{r}},$ where $p_{1}, p_{2}, \cdots p_{r}$ are primes number and $\alpha_{1},\alpha_{2}, \cdots \alpha_{r}$ positive integers. Since  $n\neq p^m$ form, so $r\ge2.$ If possible let the power graph $\mathcal{P}(C_{n})$ of cyclic group $C_{n}$ is complete multipartite for $n\neq p^m.$ Then for any two independent vertex sets $U$ and $V$ form a complete bipartite $K_{|U|,|V|}$ induced subgraph of $\mathcal{P}(C_{n}).$ Consider the independent vertex set $U$ contains the elements $\{a_{1}, a_{2}, \cdots, a_{r}\}$ such that $o(a_{i})=p_{i}.$ (since for $i\neq j,$ $o(a_{i})$ does not divide $o(a_{j})$, so neither $a_{i}^{m}=a_{j}$ nor $a_{j}^{m}=a_{i}.$ Hence $U=\{a_{1}, a_{2}, \cdots, a_{r}\}$ is an independent vertex set). And another independent vertex set $V$ contains the elements $\{b_{1}, b_{2}, \cdots, b_{r}\}$ such that $o(b_{i})=p_{i}.$ Since $o(a_{i})$ divides $o(b_{i})$, so either $a_{i}^{m}=b_{i}$ or $b_{i}^{m}=a_{i}.$ So, $ a_{i}$ is adjacent to $b_{i}.$ Since $r\ge2,$ at least one of the sets $U$ and $V$ has cardinality at least $2.$ But for $i\neq j,$ neither $a_{i}^{m}=b_{j}$ nor $b_{j}^{m}=a_{i}.$ Hence, $K_{|U|,|V|}$ is not complete bipartite induced subgraph of $\mathcal{P}(C_{n}).$ It contradicts the fact that $\mathcal{P}(C_{n})$ of cyclic group $C_{n}$ is complete multipartite for $n\neq p^m.$ Hence the power graph $\mathcal{P}(C_{n})$ of cyclic group $C_{n}$ is not complete multipartite for $n\neq p^m.$
\end{proof}

\begin{figure}[!ht]
   \begin {center}
    \includegraphics[width=2.0in]{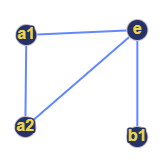}
    \caption{Induced subgraph}
   \label{fig2}
    \end{center}
\end{figure}

Again, we can improve the bound of the second largest eigenvalue of $\mathcal{P}(C_{n}),$ for the cyclic group of order $n\neq p^m.$ Here we will use the interlacing of eigenvalues and the equitable partition technique.
\begin{theorem}
    For $n\neq p^m$, the second largest eigenvalue is $\lambda_{2}(\mathcal{P}(C_{n})) \ge \mu_{2}(> 0.31),$ where $\mu_{2}$ is the second largest root of the equation $\mu^{3}-\mu^{2}-3\mu+1=0$
\end{theorem}
\begin{proof}
    Let $n=p_{1}^{\alpha_{1}}p_{2}^{\alpha_{2}}\cdots p_{r}^{\alpha_{r}},$ where $p_{1}, p_{2}, \cdots p_{r}$ are distinct primes and $\alpha_{1},\alpha_{2}, \cdots \alpha_{r}$ positive integers. Since  $n\neq p^m$ form, so $r\ge2,$ and there is odd prime $p_{i}$ such that $\phi(p_{i})\ge 2.$ Consider two elements $\{a_{1}, a_{2}\}$ of order $p_{i}$ and one element $b_{1}$ of order $p_{j},$ where $i\ne j.$ Then the elements $\{a_{1}, a_{2}, b_{1}\}$ together with the identity element $e$ form a induced subgraph of $\mathcal{P}(C_{n}),$ that induced subgraph is a graph of triangle with a pendant vertex (see Fig \ref{fig2}). We can find a quotient matrix of the induced subgraph with respect to the equitable partition $\pi=\{\{e\},\{a_{1}, a_{2}\},\{b_{1}\}\}$ is 
    \begin{eqnarray*}\displaystyle {Q_{s}}&=&{\begin{bmatrix}  {0}&{2} &{1}\\
 {1}&{1} &{0}\\
 {1}&{0} &{0}
  \end{bmatrix}}.
 \end{eqnarray*}
 The characteristic equation of the matrix $Q_{s}$ is $\mu^{3}-\mu^{2}-3\mu+1=0.$
 Then the second largest eigenvalue of the induced subgraph form by the vertices $\{e, a_{1}, a_{2},b_{1}\}$ is greater or equal to the second largest root $\mu_{2}$ of the equation $\mu^{3}-\mu^{2}-3\mu+1=0.$
 Hence, by the property of the interlacing of the eigenvalues, we get that the second largest eigenvalue is $\lambda_{2}(\mathcal{P}(C_{n})) \ge \mu_{2}.$
\end{proof}

\subsection{Dihedral Group}
Our next results provide improved spectral radius bounds for the power graphs of the dihedral groups, compared to the bounds given in \cite{srip}. The existing bounds are very trivial. Here, we present the spectral radius bounds for the power graphs  $\mathcal{P}(\mathcal{D}_{2n}),$ based on the average degree and order of the cyclic subgroups of $\mathcal{D}_{2n}.$ Also, we have found the bounds of the second largest eigenvalues of the same using the clique number of the power graph of the cyclic part.
 
\begin{theorem}\label{aD2n}
Let $\lambda_{1}(\mathcal{P}(\mathcal{D}_{2n}))$ denote the spectral radius of $\mathcal{P}(\mathcal{D}_{2n})$. For $n\ge3$, the improved spectral radius bounds are $\frac{1}{2}\left[D_{avg}+\sqrt{D^{2}_{avg}+4}\right] \le \lambda_{1}(\mathcal{P}(\mathcal{D}_{2n}))\le n$. Where $D_{avg}$ is the average degree of $\mathcal{P}(C_{n})$. 
\end{theorem}
\begin{proof} For each positive integer $n\ge 3$, the dihedral group \cite{Rose} $\mathcal{D}_{2n}=\langle a, b \rangle$ is a non-commutative group of order $2n$ where $(i) ~o(a)=n, o(b)=2$ and $(ii)~ ba=a^{-1}b=a^{n-1}b$. Since $o(a)=n, C_{n}=\langle a \rangle$ is a cyclic subgroup of $\mathcal{D}_{2n}$ of order $n$. We get the structure of $\mathcal{P}(\mathcal{D}_{2n})$ given in \cite{srip2} as a copy $\mathcal{P}(C_{n})$ and $n$ copies of the complete graph $K_{2}$ which share the identity. Let the vertices be ordered as $V_{1}=\{e, a, a^{2}, \cdots, a^{n-1}\}$ and $V_{2}=\{b, ab, a^{2}b, \cdots, a^{n-1}b \}$. So the adjacency matrix of $\mathcal{P}(\mathcal{D}_{2n})$ with rows and columns indexed by  $\{ V_{1}, V_{2}\}$ as the block matrix
    \begin{eqnarray*}\displaystyle {A(\mathcal{P}(\mathcal{D}_{2n}))}&=&{\begin{bmatrix}  {A(\mathcal{P}(C_{n}))} &{E_{n}}\\
 {E^{T}_{n}} &{O_{n}} 
 \end{bmatrix}},
 \end{eqnarray*}
 where $E_{n}=\begin{pmatrix}
     {1} &{1} &{\cdots} &{1}\\
     {0} &{0} &{\cdots} &{0}\\
     {\cdots} &{\cdots} &{\cdots} &{\cdots}\\
     {0} &{0} &{\cdots} &{0}
\end{pmatrix}$ is a square matrix of order $n$ and $A(\mathcal{P}(C_{n}))$ is given by (\ref{APC}).
Now the quotient matrix of $A(\mathcal{P}(\mathcal{D}_{2n}))$ with respect to the partition $\pi=\{ V_{1}, V_{2}\}$ is 
 \begin{eqnarray*}\displaystyle {Q_{\mathcal{D}_{2n}}}&=&{\begin{bmatrix} 
 {D_{avg}} &{1}\\
 {1} &{0} 
 \end{bmatrix}}.
 \end{eqnarray*} Where $D_{avg}$ is the average degree of $\mathcal{P}(C_{n})$. 
The spectral radius of $Q_{A}$ is $\mu_{1}({Q_{\mathcal{D}_{2n}}})=\frac{1}{2}\left[D_{avg} + \sqrt{D_{avg}^{2}+4}\right]$.
So, by the Lemma \ref{quotient}, we get $\lambda_{1}(\mathcal{P}(C_{n}))\ge \mu_{1}({Q_{\mathcal{D}_{2n}}}).$
Therefore, $\lambda_{1}(\mathcal{P}(C_{n})) \ge \frac{1}{2}\left[D_{avg} + \sqrt{D_{avg}^{2}+4}\right].$

Again, $A(\mathcal{P}(\mathcal{D}_{2n}))$ is a non-negative and irreducible matrix. So, by the Perron-Frobenius Theorem, the spectral radius $\lambda_{1}$ has multiplicity one and has an eigenvector, say $X=[x_{1}, x_{2}, \cdots, x_{2n}]^{T},$ where each $x_i$ is a positive real number.\\
Let $x_{i}=\displaystyle \max_{v_{k}\in V_{1}} x_{k}$ and $x_{j}=\displaystyle \max_{v_{k}\in V_{2}} x_{k}$. Now from the $i$th component of ${A(\mathcal{P}(\mathcal{D}_{2n}))}X=\lambda_{1}X$, we get, 
\begin{eqnarray*}
\lambda_{1} x_{i}=\displaystyle \sum_{v_{k}\in V_{1}}a_{i,k} x_{k} + \displaystyle \sum_{v_{k}\in V_{2}}a_{i,k} x_{k} \leq \displaystyle \sum_{v_{k}\in V_{1}}a_{i,k} x_{i} + \displaystyle \sum_{v_{k}\in V_{2}}a_{i,k} x_{j}=(n-1)x_{i}+nx_{j}.
\end{eqnarray*}
\begin{equation}\label{inequ1a}
    \implies (\lambda_{1}-n+1)x_{i}\leq nx_{j}
\end{equation}
Also, from the $j$th component of ${A(\mathcal{P}(\mathcal{D}_{2n}))}X=\lambda_{1} X$, we get, \begin{eqnarray*}
    \lambda_{1} x_{j}=\displaystyle \sum_{v_{k}\in V_{1}}a_{i,k} x_{k} + \displaystyle \sum_{v_{k}\in V_{2}}a_{i,k} x_{k} \leq \displaystyle \sum_{v_{k}\in V_{1}}a_{i,k} x_{i} + \displaystyle \sum_{v_{k}\in V_{2}}a_{i,k} x_{j} \leq x_{i}+ 0x_{j}
\end{eqnarray*}
\begin{equation}\label{inequ2a}
    \implies \lambda_{1}x_{j}\leq x_{i}
\end{equation}
Since $x_{i},~ x_{j}$ are positive, from inequalities (\ref{inequ1a}) \& (\ref{inequ2a}), we get
$\lambda_{1}(\lambda_{1}-n+1) \leq n.$
\begin{equation*}
\implies \lambda_{1} \le n.
\end{equation*} 
Therefore,
\begin{equation*}
    \frac{1}{2}\left[D_{avg}+\sqrt{D^{2}_{avg}+4}\right] \le \lambda_{1}(\mathcal{P}(\mathcal{D}_{2n}))\le n.
\end{equation*}
Hence, complete the proof.
\end{proof}
\textbf{Comparison:}
In the above result, although the bound is not directly comparable with the given spectral bounds in \cite{srip}, as
\begin{equation*}
\lambda_{1}(\mathcal{P}({C}_{n})) < \lambda_{1}(\mathcal{P}(\mathcal{D}_{2n})) \le \lambda_{1}(\mathcal{P}({C}_{n})) + \sqrt{n},
\end{equation*}
our bounds are much better than the above bounds. We can observe this through the following two examples:
  
\begin{example}
    \emph{Let us consider the group $\mathbb{G}=\mathbb{Z}_{6}$. Then the spectral radius of $A(\mathcal{P}(\mathbb{Z}_{6}))$ is the largest root of the equation, $x^{3}-3x^{2}-7x+3=0$ i.e., $\lambda_{1}(\mathcal{P}(\mathbb{Z}_{6})) \approx 4.42788$. Now the spectral radius bound of the power graph of the group   $\mathbb{G}=\mathcal{D}_{12}$ is $4.42788< \lambda_{1}(\mathcal{P}(\mathcal{D}_{12}))\le 4.42788+\sqrt{6}$, according to the bound given in \cite{srip}. i.e., 
   \begin{eqnarray*}
       4.42788<\lambda_{1}(\mathcal{P}(\mathcal{D}_{12}))\le 6.87737.
   \end{eqnarray*} 
    For the graph $\mathcal{P}(\mathbb{Z}_{6})$, we get $D_{avg}=\frac{13}{3}$ (from Fig\ref{fig1}) and $\frac{1}{2}\left[D_{avg}+\sqrt{D^{2}_{avg}+4}\right]=\frac{1}{2}\left[\frac{13}{3}+\sqrt{(\frac{13}{3})^{2}+4}\right]\approx 4.55297$. So, by our result in Theorem \ref{aD2n}, the spectral radius bound of the power graph of the group   $\mathbb{G}=\mathcal{D}_{12}$ is 
   \begin{eqnarray*}
       4.55297\le \lambda_{1}(\mathcal{P}(\mathcal{D}_{12}))\le 6.
   \end{eqnarray*} }
\end{example}
\begin{example}
  \emph{ For $n=p^m$, for any prime $p$ and any natural number $m$, the spectral radius bound of the power graph of the group   $\mathbb{G}=\mathcal{D}_{2p^m}$ is $\lambda_{1}(\mathcal{P}({C}_{p^m}))< \lambda_{1}(\mathcal{P}(\mathcal{D}_{2p^m})) \le \lambda_{1}(\mathcal{P}({C}_{p^m})) + \sqrt{p^m}$, according to the bound given in \cite{srip}. i.e.
   \begin{eqnarray*}
       p^{m}-1< \lambda_{1}(\mathcal{P}(\mathcal{D}_{2p^m}))\le p^{m} + \sqrt{p^m}-1.
   \end{eqnarray*} 
  But our result provides a comparably better bound, which is
  \begin{eqnarray*}
      \frac{1}{2}\left[p^{m}-1+\sqrt{(p^{m}-1)^{2}+4}\right] \le \lambda_{1}(\mathcal{P}(\mathcal{D}_{2p^m}))\le p^{m}.
  \end{eqnarray*}}
\end{example}
Now we find the lower bound of the second largest eigenvalue of $\mathcal{P}(\mathcal{D}_{2n}).$
\begin{theorem}
\label{2nd_aD2n}
 For $n\ge3$, the second largest eigenvalue is $ \lambda_{2}(\mathcal{P}(\mathcal{D}_{2n}))\ge \mu_{2},$ where $\mu_{2}$ is the second largest root of the equation, $\mu^{3}-(\omega-2)\mu^{2}-(n+\omega-1)\mu+(\omega-2)n=0$ and $\omega=\omega(\mathcal{P}(\mathcal{D}_{2n}))$ is the clique number of $\mathcal{P}(\mathcal{D}_{2n}).$
\end{theorem}
\begin{proof}
For $n\ge3, C_{n}=\langle a \rangle$ is a cyclic subgroup of $\mathcal{D}_{2n}=\langle a, b \rangle$ of order $n$. We know the structure of $\mathcal{P}(\mathcal{D}_{2n})$ as a copy $\mathcal{P}(C_{n})$ and $n$ copies of the complete graph $K_{2}$ that share the identity. Let $\omega=\omega(\mathcal{P}(C_{n}))=\omega(\mathcal{P}(\mathcal{D}_{2n}))$ be the clique number of $\mathcal{P}(C_{n})$. This clique, together with $n$ copies of the complete graph $K_{2}$ which share the identity, form an induced subgraph of the graph $\mathcal{P}(\mathcal{D}_{2n}).$ Let the vertices of this induced subgraph is $U=U_{1}\cup U_{2}\cup U_{3},$ where $U_{1}=\{e\}, ~U_{2}=\{a_{2}, a_{3}, \cdots a_{\omega}\}, ~U_{3}=\{b, ab, a^{2}b, \cdots, a^{n-1}b \}.$
We can find a quotient matrix of the induced subgraph with respect to the equitable partition $\pi=\{U_{1}, U_{2}, U_{3}\}$ is 
    \begin{eqnarray*}\displaystyle {Q_{s}}&=&{\begin{bmatrix}  {0} &{\omega-1} &{n}\\
 {1}&{\omega-2} &{0}\\
 {1}&{0} &{0}
 \end{bmatrix}}.
 \end{eqnarray*}
 The characteristic equation of the matrix $Q_{s}$ is $\mu^{3}-(\omega-2)\mu^{2}-(n+\omega-1)\mu+(\omega-2)n=0.$  Then the second largest eigenvalue of the induced subgraph form by the vertices $U=U_{1}\cup U_{2}\cup U_{3}$ is greater or equal to the second largest root $\mu_{2}$ of the equation $\mu^{3}-(\omega-2)\mu^{2}-(n+\omega-1)\mu+(\omega-2)n=0.$
 Hence, by the property of the interlacing of the eigenvalues, we get that the second largest eigenvalue is $\lambda_{2}(\mathcal{P}(\mathcal{D}_{2n})) \ge \mu_{2}.$
\end{proof}
\begin{corollary}
    The graph $\mathcal{P}(\mathcal{D}_{2n})$ has at least $n$ number of non negative eigenvalues.
\end{corollary}
\begin{proof}
    From the structure of $\mathcal{P}(\mathcal{D}_{2n})$ we get $n$ copies of the complete graph $K_{2}$ which share the identity. So there is an independent vertex set of order $n.$ So there is a zero matrix $O_{n \times n}$ of order $n,$ is a principal submatrix of $A(\mathcal{P}(\mathcal{D}_{2n})).$ Hence by Cauchy interlacing, $\lambda_{1}(\mathcal{P}(\mathcal{D}_{2n}))\ge \lambda_{2}(\mathcal{P}(\mathcal{D}_{2n}))\ge \cdots \ge \lambda_{n}(\mathcal{P}(\mathcal{D}_{2n})) \ge 0.$
\end{proof}


\subsection{Dicyclic Group}
\begin{figure}[!ht]
   \begin {center}
    \includegraphics[width=3.0in]{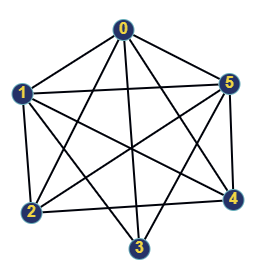}
    \caption{Power graph of $\mathbb{Z}_{6}$}
   \label{fig1}
    \end{center}
\end{figure}

The next results provide improved spectral radius bounds for the power graphs of the dicyclic groups, compared to the bounds given in \cite{srip}. The earlier bounds are very trivial. Here, we present the spectral radius bounds for the power graph $\mathcal{P}(\mathcal{Q}_{4n})$, based on the average degree and order of the cyclic subgroup $\mathcal{Q}_{4n}$. Also, we have included the bounds of the second largest eigenvalues of the same.
\begin{theorem}\label{aQ4n}
Let $\lambda_{1}(\mathcal{P}(\mathcal{Q}_{4n}))$ denote the spectral radius of $\mathcal{P}(\mathcal{Q}_{4n})$. For any integer $n\ge2$, the improved bounds on spectral radius are $\frac{1}{2}\left[(D_{avg}+1)+\sqrt{(D_{avg}-1)^{2}+16}\right] \le \lambda_{1}(\mathcal{P}(\mathcal{Q}_{4n}))\le 2n+1 $. Where $D_{avg}$ is the average degree of $\mathcal{P}(C_{2n})$.
\end{theorem}
\begin{proof}  For each positive integer $n\ge 2$, the dicyclic group \cite{Rose} $\mathcal{Q}_{4n}=\langle a, b \rangle$ is a non-commutative group of order $4n$ where $(i) a^{2n}=e, a^{n}=b^{2}$ and $(ii) ab=ba^{-1}=ba^{2n-1}$. So for $1\le i \le n$, $A_{i}=\{ e, a^{i-1}b, a^{i}, a^{n+i-1}b \}$ and $C_{2n}=\langle a \rangle$ are the only maximal cyclic subgroups of $\mathcal{Q}_{4n}$. Then we get the structure of $\mathcal{P}(\mathcal{Q}_{4n})$ given in \cite{srip2}, as a copy of $\mathcal{P}(C_{2n})$ and $n$ copies of the complete graph $K_{4}$ (each of which is $\mathcal{P}(A_{i}$)) sharing the identity $e$ and the unique involution $a^{n}$. Let the vertices be ordered as $V_{1}=\{e, a^{n}, a, a^{2}, \cdots, a^{n-1}, a^{n+1}, \cdots, a^{2n-1}\}$ and $V_{2}=\{b, ab, a^{2}b, \cdots, a^{2n-1}b \}$. We get the adjacency matrix of $\mathcal{P}(C_{2n})$ in terms of block matrix as 
    \begin{eqnarray*}\displaystyle {A(\mathcal{P}(\mathcal{Q}_{4n}))}&=&{\begin{bmatrix}  {A(\mathcal{P}(C_{2n}))} &{F_{2n}}\\
 {F^{T}_{2n}} &{P_{2n}} 
 \end{bmatrix}},
 \end{eqnarray*}
 where $F_{2n}=\begin{pmatrix}
     {1} &{1} &{\cdots} &{1}\\
     {1} &{1} &{\cdots} &{1}\\
     {0} &{0} &{\cdots} &{0}\\
     {\cdots} &{\cdots} &{\cdots} &{\cdots}\\
     {0} &{0} &{\cdots} &{0}
 \end{pmatrix}$ is a square matrix of order $2n$ and $P_{2n}$ is the block matrix $\begin{bmatrix}
     {O_{n}} &{I_{n}}\\
     {I_{n}} &{O_{n}}
 \end{bmatrix}.$
 Now the quotient matrix of $A(\mathcal{P}(\mathcal{Q}_{4n}))$ with respect to the partition $\pi=\{ V_{1}, V_{2}\}$ is 
 \begin{eqnarray*}\displaystyle {Q_{\mathcal{Q}_{4n}}}&=&{\begin{bmatrix} 
 {D_{avg}} &{2}\\
 {2} &{1} 
 \end{bmatrix}}.
 \end{eqnarray*} Where $D_{avg}$ is the average degree of $\mathcal{P}(C_{2n})$.

The spectral radius of $Q_{\mathcal{Q}_{4n}}$ is $\mu_{1}({Q_{\mathcal{Q}_{4n}}})=\frac{1}{2}\left[(D_{avg}+1) + \sqrt{(D_{avg}-1)^{2}+16}\right].$
So, by the Lemma \ref{quotient}, we get $\lambda_{1}(\mathcal{Q}_{4n})\ge \mu_{1}({Q_{\mathcal{Q}_{4n}}}).$
Therefore, $\lambda_{1}(\mathcal{Q}_{4n}) \ge \frac{1}{2}\left[(D_{avg}+1) + \sqrt{(D_{avg}-1)^{2}+16}\right].$

Again, as $A(\mathcal{P}(\mathcal{Q}_{4n}))$ is non-negative and irreducible, by the Perron-Frobenius Theorem, the spectral radius $\lambda_{1}$ has multiplicity one and has an eigenvector, say $X=[x_{1}, x_{2}, \cdots, x_{4n}]^{T},$ where each $x_i$ is a positive real number.\\
Let $x_{i}=\displaystyle \max_{v_{k}\in V_{1}} x_{k}$ and $x_{j}=\displaystyle \max_{v_{k}\in V_{2}} x_{k}$. Now from the $i$th component of ${A(\mathcal{P}(\mathcal{Q}_{4n}))}X=\lambda_{1}X$, we get, 
\begin{eqnarray*}
\lambda_{1} x_{i}=\displaystyle \sum_{v_{k}\in V_{1}}a_{i,k} x_{k} + \displaystyle \sum_{v_{k}\in V_{2}}a_{i,k} x_{k} \leq \displaystyle \sum_{v_{k}\in V_{1}}a_{i,k} x_{i} + \displaystyle \sum_{v_{k}\in V_{2}}a_{i,k} x_{j}=(2n-1)x_{i}+2nx_{j}.
\end{eqnarray*}
\begin{equation}\label{inequ1q}
    \implies (\lambda_{1}-2n+1)x_{i}\leq 2nx_{j}
\end{equation}
Again from the $j$th component of ${A(\mathcal{P}(\mathcal{Q}_{4n}))}X=\lambda_{1} X$, we get, \begin{eqnarray*}
    \lambda_{1} x_{j}=\displaystyle \sum_{v_{k}\in V_{1}}a_{i,k} x_{k} + \displaystyle \sum_{v_{k}\in V_{2}}a_{i,k} x_{k} \leq \displaystyle \sum_{v_{k}\in V_{1}}a_{i,k} x_{i} + \displaystyle \sum_{v_{k}\in V_{2}}a_{i,k} x_{j} \leq 2x_{i}+ x_{j}
\end{eqnarray*}
\begin{equation}\label{inequ2q}
    \implies (\lambda_{1}-1)x_{j}\leq 2x_{i}
\end{equation}
Since $x_{i},~ x_{j}$ are positive, from inequality (\ref{inequ1q}) \& (\ref{inequ2q}), we get
$(\lambda_{1}-1)(\lambda_{1}-2n+1) \leq 4n.$

\begin{equation*}
\implies \lambda_{1} \le 2n+1.
\end{equation*} Therefore,
\begin{equation*}
\frac{1}{2}\left[(D_{avg}+1)+\sqrt{(D_{avg}-1)^{2}+16}\right] \le \lambda_{1}(\mathcal{P}(\mathcal{Q}_{4n}))\le 2n+1.
\end{equation*}
Hence, complete the proof.
\end{proof}
\textbf{Comparison:} In the same way as earlier, the bounds are not directly comparable with the given spectral bounds in \cite{srip}, as
\begin{equation*}
\lambda_{1}(\mathcal{P}(C_{2n})) < \lambda_{1}(\mathcal{P}(\mathcal{Q}_{4n})) \le \lambda_{1}(\mathcal{P}(C_{2n})) + 2\sqrt{n}.
\end{equation*}
However, our bounds are significantly better than the above bounds, which we can observe with the following two examples:

\begin{example}
 \emph{For $n=3$, the spectral radius bound of the power graph of the group   $\mathbb{G}=\mathcal{Q}_{12}$ is $\lambda_{1}(\mathcal{P}(\mathbb{Z}_{6}))< \lambda_{1}(\mathcal{P}(\mathcal{Q}_{12})) \le \lambda_{1}(\mathcal{P}(\mathbb{Z}_{6})) + 2\sqrt{3}$, according to the bound given in \cite{srip}. i.e., 
   \begin{eqnarray*}
       4.42788<\lambda_{1}(\mathcal{P}(\mathcal{Q}_{12}))\le 7.89198.
   \end{eqnarray*} 
For the graph $\mathcal{P}(\mathbb{Z}_{6})$, we get $D_{avg}=\frac{13}{3}$ (from Fig \ref{fig1}) and $\frac{1}{2}\left[(D_{avg}+1)+\sqrt{(D_{avg}-1)^{2}+16}\right] \\= \frac{1}{2}\left[\frac{16}{3}+\sqrt{(\frac{10}{3})^{2}+16}\right] \approx 5.27008$. So, by our result in Theorem \ref{aQ4n}, the spectral radius bound of the power graph of the group   $\mathbb{G}=\mathcal{Q}_{12}$ is 
   \begin{eqnarray*}
       5.27008\le \lambda_{1}(\mathcal{P}(\mathcal{Q}_{12}))\le 7.
   \end{eqnarray*} }
\end{example}
\begin{example}
   \emph{ For $n=2^m$, for any natural number $m\ge2$, the spectral radius bound of the power graph of the dicyclic group   $\mathbb{G}=\mathcal{Q}_{2^{m+2}}$ is $\lambda_{1}(\mathcal{P}({C}_{2^{m+1}}))< \lambda_{1}(\mathcal{P}(\mathcal{Q}_{2^{m+2}})) \le \lambda_{1}(\mathcal{P}({C}_{2^{m+1}})) + 2\sqrt{2^m}$, according to the bound given in \cite{srip}. i.e.,
   \begin{eqnarray*}
       2^{m+1}-1< \lambda_{1}(\mathcal{P}(\mathcal{Q}_{2^{m+2}}))\le 2^{m+1} + 2\sqrt{2^m}-1.
   \end{eqnarray*} 
  But our result provides a comparably better bound, which is
  \begin{eqnarray*}
      \frac{1}{2}\left[2^{m+1}+\sqrt{(2^{m+1}-2)^{2}+16}\right] \le \lambda_{1}(\mathcal{P}(\mathcal{Q}_{2^{m+2}}))\le 2^{m+1}+1\\
      \implies \left[2^{m}+\sqrt{(2^{m}-1)^{2}+4}\right] \le \lambda_{1}(\mathcal{P}(\mathcal{Q}_{2^{m+2}}))\le 2^{m+1}+1.
  \end{eqnarray*}}
  \end{example}
  Next, we give the second largest eigenvalue lower bound of $\mathcal{P}(\mathcal{Q}_{4n}).$
\begin{theorem}\label{2ndaQ4n}
 For $n\ge2$, the second largest eigenvalue is $ \lambda_{2}(\mathcal{P}(\mathcal{Q}_{4n}))\ge \mu_{2},$ where $\mu_{2}$ is the second largest root of the equation, $\mu^{3}-(k+1)\mu^{2}-(4n+1)\mu+(4kn-4n+k+1)=0$ and $k=\phi(2n)$ be the number of generators of $C_{2n}.$
\end{theorem}
\begin{proof} 
For $n\ge2, C_{2n}=\langle a \rangle$ is a cyclic subgroup of $\mathcal{Q}_{4n}=\langle a, b \rangle$ of order $2n$. We know the structure of $\mathcal{P}(\mathcal{Q}_{4n})$ given in \cite{srip2}, as a copy of $\mathcal{P}(C_{2n})$ and $n$ copies of the complete graph $K_{4}$  sharing the identity $e$ and the unique involution $a^{n}$. Let $k=\phi(2n)$ be the number of generators of $C_{2n}.$ Then the vertices $U=U_{1}\cup U_{2}\cup U_{3},$ where $U_{1}=\{e, a^{n}\}, ~U_{2}=\{a_{1}, a_{2}, \cdots a_{k}\}, ~U_{3}=\{b, ab, a^{2}b, \cdots, a^{2n-1}b \},$ form an induced subgraph of $\mathcal{P}(\mathcal{Q}_{4n}).$ We can find a quotient matrix of the induced subgraph with respect to the equitable partition $\pi=\{U_{1}, U_{2}, U_{3}\}$ is 
    \begin{eqnarray*}\displaystyle {Q_{s}}&=&{\begin{bmatrix}  {1} &{k} &{2n}\\
 {2}&{k-1} &{0}\\
 {2}&{0} &{1}
 \end{bmatrix}}.
 \end{eqnarray*}
 The characteristic equation of the matrix $Q_{s}$ is $\mu^{3}-(k+1)\mu^{2}-(4n+1)\mu+(4kn-4n+k+1)=0.$ Then the second largest eigenvalue of the induced subgraph is greater or equal to the second largest root $\mu_{2}$ of the equation $\mu^{3}-(k+1)\mu^{2}-(4n+1)\mu+(4kn-4n+k+1)=0.$
 Hence, by the property of the interlacing of the eigenvalues, we get that the second largest eigenvalue is $\lambda_{2}(\mathcal{P}(\mathcal{Q}_{4n})) \ge \mu_{2}.$
\end{proof}
\begin{corollary}
    For $n=2^{m}$, the second largest eigenvalue is $ \lambda_{2}(\mathcal{P}(\mathcal{Q}_{4n}))\ge \mu_{2},$ where $\mu_{2}$ is the second largest root of the equation, $\mu^{3}-(\omega+1)\mu^{2}-(2\omega +1)\mu+(2\omega^{2}-5\omega -1)=0$ and $\omega=2n$ be the clique number of $\mathcal{P}(\mathcal{Q}_{4n}).$
\end{corollary}
\begin{proof}
    For $n=2^{m}, \mathcal{P}(C_{2n})=K_{2n}.$ So, the clique number of $\mathcal{P}(\mathcal{Q}_{4n}) $ is $\omega=2n.$ Let $U_{1}=\{e, a^{n}\}, ~U_{2}=\{a, a^{2}, \cdots a^{n-1}, a^{n+1}, \cdots, a^{2n-1} \}, ~U_{3}=\{b, ab, a^{2}b, \cdots, a^{2n-1}b \},$ are the vertices of the graph $\mathcal{P}(\mathcal{Q}_{4n}).$ The quotient matrix with respect to the equitable partition $\pi=\{U_{1}, U_{2}, U_{3}\}$ is 
     \begin{eqnarray*}\displaystyle {Q_{s}}&=&{\begin{bmatrix}  {1} &{\omega-2} &{\omega}\\
 {2}&{\omega-3} &{0}\\
 {2}&{0} &{1}
 \end{bmatrix}}.
 \end{eqnarray*}
 The characteristic equation of the matrix $Q_{s}$ is $\mu^{3}-(\omega+1)\mu^{2}-(2\omega +1)\mu+(2\omega^{2}-5\omega -1)=0.$ Then the second largest eigenvalue of the graph is greater or equal to the second largest root $\mu_{2}$ of the equation $\mu^{3}-(\omega+1)\mu^{2}-(2\omega +1)\mu+(2\omega^{2}-5\omega -1)=0.$
 Hence, by the property of the interlacing of the eigenvalues, we get that the second largest eigenvalue is $\lambda_{2}(\mathcal{P}(\mathcal{Q}_{4n})) \ge \mu_{2}.$
\end{proof}

Next, we find the spectral radius of the family of cyclic groups of order $pq$, where $p$ and $q$ are primes, expressed in terms of the largest root of a cubic equation.

\begin{theorem}\label{Zpq}
    Let $\mathbb{G}=\mathbb{Z}_{pq}$, where $p$ and $q$ be primes; then the spectral radius of $A(\mathcal{P}(\mathbb{Z}_{pq}))$ is the largest root of the equation $x^{3}-(pq-3)x^{2}-(pq+p+q-4)x+(p^{2}q^{2}-2p^{2}q-2pq^{2}+p^{2}+5pq+q^{2}-4p-4q+4)=0$.
\end{theorem}
\begin{proof}
     Let us take a partition of the vertex set as $\{ V_{1}, V_{2}, V_{3}\}$ such that $V_{1}$ contains all generators along with the identity element of $\mathbb{Z}_{pq}$ (so, $|V_{1}|= \phi(pq) + 1 = (p-1)(q-1)+1$), $V_{2}=\{ p, 2p, \cdots, (q-1)p\}$ and $V_{3}=\{ q, 2q, \cdots, (p-1)q\}$. Then for all \( a, b \) in the set $V_{2}\cup\{0\}$, the sum $( a + b)\mod pq$ is also in $V_{2}\cup\{0\}$. Therefore, this set is a subgroup $\mathbb{Z}_{pq}$ of order $q$. Then its power graph $\mathcal{P}({V_{2}\cup\{0\}})$ is a complete graph of order $q$. With a similar argument, we have the power graph, which $\mathcal{P}({V_{3}\cup\{0\}})$ is a complete graph of order $p$. Then the subgraph $\mathcal{P}(V_{2})$ of $\mathcal{P}({V_{2}\cup\{0\}})$ is a complete graph of order $q-1$ and the subgraph $\mathcal{P}(V_{3})$ of $\mathcal{P}({V_{3}\cup\{0\}})$ is a complete graph of order $p-1$. Now we prove any vertex in $\mathcal{P}(V_{2})$ is not adjacent to any vertex in $\mathcal{P}(V_{3})$. Let $a\in \{ 1, 2, \cdots, (q-1)\}$ and $ b\in \{ 1, 2, \cdots, (p-1)\}$ be two arbitrary elements. Then $ap$ is a vertex in $\mathcal{P}(V_{2})$ and $bq$ is a vertex in $\mathcal{P}(V_{3})$. If possible, let $ap$ and $bq$ be adjacent in $\mathcal{P}(\mathbb{Z}_{pq})$. Then there exist positive integers $m$ and $n$ such that 
  \begin{eqnarray*}
      ap=mbq ~~\mbox{or}~~ bq=nap.
  \end{eqnarray*}
  Now if $ap=mbq$, then $ap=mbq+kpq$ for some $k\in \mathbb{Z}$. Then it shows that the prime $q$ divides $ap$, which is impossible as $q \nmid p$ and $q \nmid a$. Therefore, $ap\neq mbq$. Similarly, we can say that $bq \neq nap$. Therefore, $ap$ and $bq$ are not adjacent. So, we can write the adjacency matrix of the power graph $\mathcal{P}(\mathbb{Z}_{pq})$ with rows and columns indexed by $\{ V_{1}, V_{2}, V_{3}\}$:
   \begin{eqnarray*} \displaystyle {A(\mathcal{P}(\mathbb{Z}_{pq}))}&=&{\begin{bmatrix} 
    {J_{\phi(pq)+1}-I_{\phi(pq)+1}} & {J_{(\phi(pq)+1)\times (q-1)}} & {J_{(\phi(pq)+1)\times (p-1)}}\\
 {J_{(q-1)\times (\phi(pq)+1)}} & {J_{(q-1)}-I_{(q-1)}} & {O_{(q-1)\times (p-1)}}\\
 {J_{(p-1) \times (\phi(pq)+1)}} & {O_{(p-1)\times (q-1)}} & {J_{(p-1)}-I_{(p-1)}}
 \end{bmatrix}}.
 \end{eqnarray*}
 
 Then the quotient matrix of $A(\mathcal{P}(\mathbb{Z}_{pq}))$ with respect to the partition $\pi=\{ V_{1}, V_{2}, V_{3}\}$ is 
 \begin{eqnarray*}\displaystyle {Q_{A}^{\mathcal{P}(\mathbb{Z}_{pq})}}&=&{\begin{bmatrix} 
 {(p-1)(q-1)} &{q-1} &{p-1}\\
 {(p-1)(q-1)+1} &{q-2} &{0}\\
 {(p-1)(q-1)+1} &{0} &{p-2}
 \end{bmatrix}}.
\end{eqnarray*} It can be observed that the partition $\pi=\{ V_{1}, V_{2}, V_{3}\}$ is an equitable partition. Therefore, the spectral radius of the matrix $A(\mathcal{P}(\mathbb{Z}_{pq}))$ is nothing but the spectral radius of the matrix $Q_{A}^{\mathcal{P}(\mathbb{Z}_{pq})}$. The (Perron-Frobenius) Theorem \ref{perron} and Theorem \ref{spec} ensure that the spectral radius of the matrix  $A(\mathcal{P}(\mathbb{Z}_{pq}))$ is the largest root of the characteristic polynomial of the matrix $Q_{A}^{\mathcal{P}(\mathbb{Z}_{pq})}$ and the required characteristic equation is det$(xI-Q_{A}^{\mathcal{P}(\mathbb{Z}_{pq})})=0$ i.e., $x^{3}-(pq-3)x^{2}-(pq+p+q-4)x+(p^{2}q^{2}-2p^{2}q-2pq^{2}+p^{2}+5pq+q^{2}-4p-4q+4)=0$.
\end{proof}



\section{Bounds on the distance spectral radius }

In a power graph of a group, the identity element (vertex) is adjacent to every other element (vertex), so the power graph of a group is always connected and has a diameter at most $2$. Thus, we can consider the distance matrix of the power graph of a finite group. In this section, we find similar bounds for the distance spectral radius of the power graph of certain finite groups, applying the same technique. The following result provides the bounds for the distance spectral radius of the power graph of a finite cyclic group.
\begin{theorem}\label{dCn}
Let $\rho_{1}(\mathcal{P}(C_{n}))$ denote the spectral radius of $\mathcal{P}(C_{n})$. For any integer $n\ge2$, the spectral radius $\rho_{1}(\mathcal{P}(C_{n}))$ satisfies: $$ \rho_{1}(\mathcal{P}(C_{n})) \ge \frac{1}{2}\left[ (Tr_{avg}-1) + \sqrt{(Tr_{avg}-2l+1)^{2}+4l(n-l)}\right] $$ and $$\rho_{1}(\mathcal{P}(C_{n}))\le \frac{1}{2}\left[(Tr_{max}-1) + \sqrt{(Tr_{max}-2l+1)^{2}+4l(n-l)}\right].$$ Where $l=\phi(n)+1$, $Tr_{avg}$, and $Tr_{max}$ are the average transmissions and maximum transmissions of the non-identity non-generator elements of $\mathcal{P}(C_{n})$, respectively. Moreover, equality holds in both bounds if and only if $n=p^{m}$, for any prime number $p$ and any positive integer $m$.
\end{theorem}
\begin{proof}
Let $V_{1}$, and $V_{2}$ be two partition sets of the vertices of the power graph of the cyclic group $C_{n}$. Sets $V_{1}$ contain the identity element and all generators of $C_{n}$, and $V_{2}$ contains the rest of the vertices. Then $ |V_{1}|= \phi(n)+1=l$ (say), where $\phi(n)$ is Euler's $\phi$ function. Here we are ordering the vertices of $V_{1}$ as $\{v_{1}=e, v_{2}, \cdots, v_{l}\}$, and $V_{2}$ as $\{ v_{l+1}, v_{l+2}, \cdots, v_{n}\}$. Now the adjacency matrix of the power graph $\mathcal{P}(C_{n})$ of the cyclic group $C_{n}$ with rows and columns indexed by $\{ V_{1}, V_{2}\}$ is 
     \begin{eqnarray} \label{DPC}\displaystyle {D(\mathcal{P}(C_{n}))}&=&{\begin{bmatrix} 
    {J_{l\times l}-I_{l\times l}} & {J_{l\times (n-l)}}\\
 {J_{(n-l)\times l}} &{D(\mathcal{P}(V_{2}))_{n-l}} 
 \end{bmatrix}}.
 \end{eqnarray}
 Where $J$ is a matrix with all entries $1$ with suitable order, $\textbf{1}$ is the vector with all entries $1$ with suitable order,  $I$ is the identity matrix, and $D(\mathcal{P}(V_{2}))=(d_{i,j})$ is the distance matrix of the power graph induced by the vertices of $V_{2}.$ Now the quotient matrix of $D(\mathcal{P}(C_{n}))$ matrix with respect to the partition $\pi=\{ V_{1}, V_{2}\}$ is 
 \begin{eqnarray*}\displaystyle {Q_{D}}&=&{\begin{bmatrix} 
 {l-1} &{n-l}\\
 {l} &{Tr_{avg}-l} 
 \end{bmatrix}}.
 \end{eqnarray*}
 
 Then the spectral radius of $Q_{D}$ is $\mu_{1}({Q_{D}})=\frac{1}{2}\left[(Tr_{avg}-1) + \sqrt{(Tr_{avg}-2l+1)^{2}+4l(n-l)}\right].$
So, by Lemma \ref{quotient}, we get $\rho_{1}(\mathcal{P}(C_{n}))\ge \mu_{1}({Q_{D}}).$\\
Therefore, $\rho_{1}(\mathcal{P}(C_{n})) \ge \frac{1}{2}\left[(Tr_{avg}-1) + \sqrt{(Tr_{avg}-2l+1)^{2}+4l(n-l)}\right].$ 

Again, ${D(\mathcal{P}(C_{n}))}$ is a non-negative and irreducible matrix. So, by the Perron-Frobenius Theorem, the spectral radius $\rho_{1}$ is simple and has a positive eigenvector, say $X=[x_{1}, x_{2}, \cdots, x_{n}]^{T}.$\\
Let $x_{i}=\displaystyle \max_{v_{k}\in V_{1}} x_{k}$ and $x_{j}=\displaystyle \max_{v_{k}\in V_{2}} x_{k}$. Then, considering the $i$th component of ${D(\mathcal{P}(C_{n}))}X=\rho_{1}X$, we get, 
\begin{eqnarray*}
\rho_{1} x_{i}=\displaystyle \sum_{v_{k}\in V_{1}}d_{i,k} x_{k} + \displaystyle \sum_{v_{k}\in V_{2}}d_{i,k} x_{k} \leq \displaystyle \sum_{v_{k}\in V_{1}}d_{i,k} x_{i} + \displaystyle \sum_{v_{k}\in V_{2}}d_{i,k} x_{j}=(l-1)x_{i}+(n-l)x_{j}.
\end{eqnarray*}
\begin{equation}\label{inequ1}
    \implies (\rho_{1}-l+1)x_{i}\leq (n-l)x_{j}
\end{equation}
Also, considering the $j$th component of ${D(\mathcal{P}(C_{n}))}X=\rho_{1} X$, we get, \begin{eqnarray*}
    \rho_{1} x_{j}=\displaystyle \sum_{v_{k}\in V_{1}}d_{i,k} x_{k} + \displaystyle \sum_{v_{k}\in V_{2}}d_{i,k} x_{k} \leq \displaystyle \sum_{v_{k}\in V_{1}}d_{i,k} x_{i} + \displaystyle \sum_{v_{k}\in V_{2}}d_{i,k} x_{j} \leq lx_{i}+(Tr_{max}-l)x_{j}
\end{eqnarray*}
\begin{equation}\label{inequ2}
    \implies (\rho_{1}-Tr_{max}+l)x_{j}\leq lx_{i}
\end{equation}
Since $x_{i},~ x_{j}$ are positive, from inequality (\ref{inequ1}) \& (\ref{inequ2}), we get
$(\rho_{1}-l+1)(\rho_{1}-Tr_{max}+l) \leq (n-l)l.$
\begin{equation*}
    \implies \rho_{1}^{2}-\left[(l-1)+(Tr_{max}-l)\right]\rho_{1}+(l-1)(Tr_{max}-l)-(n-l)l\leq 0
\end{equation*}
As $\rho_{1} >0$
\begin{equation*}
    \rho_{1}(\mathcal{P}(C_{n})) \leq \frac{1}{2}\left[ (Tr_{max}-1) + \sqrt{(Tr_{max}-2l+1)^{2}+4l(n-l)}\right].
\end{equation*}
The both side equality holds of $\rho_{1}(\mathcal{P}(C_{n}))$ if and only if $Tr_{avg}=Tr_{max}$, this can happen only when $\mathcal{P}(C_{n})=K_{n}$. So, the equality holds for the bounds if and only if $n=p^m$, for any prime number $p$ and any positive integer $m$.
\end{proof}

Now we provide the distance spectral radius bounds of the power graph of the dihedral group $\mathcal{D}_{2n}$ depending on the average transmission and maximum transmission of the cyclic subgroup of $\mathcal{D}_{2n}$.
\begin{theorem}\label{dD2n}
For any integer $n\ge 3$, the distance spectral radius $\rho_{1}(\mathcal{P}(\mathcal{D}_{2n}))$ of $\mathcal{P}(\mathcal{D}_{2n})$ satisfies the following: $$ \rho_{1}(\mathcal{P}(\mathcal{D}_{2n}))\ge \frac{1}{2}\left[(T_{avg}+2n-2) + \sqrt{(T_{avg}-2n+2)^{2}+4(2n-1)^{2}}\right] $$ and $$ \rho_{1}(\mathcal{P}(\mathcal{D}_{2n}))\le \frac{1}{2}\left[(T_{max}+2n-2) + \sqrt{(T_{max}-2n+2)^{2}+8n(2n-1)}\right]$$ where $T_{avg}$ and $T_{max}$ are the average and maximum transmissions of the graph $\mathcal{P}(C_{n})$ respectively.
\end{theorem}
\begin{proof} We are following the same vertex labeling and ordering as in Theorem \ref{aD2n}. So the distance matrix of $\mathcal{P}(\mathcal{D}_{2n})$ with rows and columns indexed by  $\{ V_{1}, V_{2}\}$ as the block matrix
    \begin{eqnarray*}\displaystyle {D(\mathcal{P}(\mathcal{D}_{2n}))}&=&{\begin{bmatrix}  {D(\mathcal{P}(C_{n}))} &{B_{n}}\\
 {B^{T}_{n}} &{2J_{n}-2I_{n}} 
 \end{bmatrix}},
 \end{eqnarray*} where $B_{n}$ is the block matrix $\begin{bmatrix}
     {\textbf{1}^{T}} \\
     {2J_{(n-1)\times n}} 
      \end{bmatrix}$ and $D(\mathcal{P}(C_{n}))$ is given by (\ref{DPC}). Then the quotient matrix of $D(\mathcal{P}(\mathcal{D}_{2n}))$ with respect to the partition $\pi=\{ V_{1}, V_{2}\}$ is 
\begin{eqnarray*}\displaystyle {Q_{D}}&=&{\begin{bmatrix} 
 {T_{avg}} &{2n-1}\\
 {2n-1} &{2n-2} 
 \end{bmatrix}};
\end{eqnarray*} where $T_{avg}$ is the average transmission of $\mathcal{P}(C_{n})$. The spectral radius of $Q_{D}$ is \\ $\mu_{1}({Q_{D}})=\frac{1}{2}\left[(T_{avg} + 2n-2)+ \sqrt{(T_{avg}-2n+2)^{2}+4(2n-1)^{2}}\right]$.
So, by the Lemma \ref{quotient}, we get $\rho_{1}(\mathcal{P}(\mathcal{D}_{2n}))\ge \mu_{1}({Q_{D}}).$
Therefore, $\rho_{1}(\mathcal{P}(\mathcal{D}_{2n}))\ge \frac{1}{2}\left[(T_{avg} + 2n - 2) + \sqrt{(T_{avg}-2n +2)^{2}+4(2n-1)^{2}}\right].$

Again, $D(\mathcal{P}(\mathcal{D}_{2n}))$ is a non-negative and irreducible matrix. So, by the Perron-Frobenius Theorem, the spectral radius $\rho_{1}$ is simple and has a positive eigenvector, say $X=[x_{1}, x_{2}, \cdots, x_{2n}]^{T}.$\\
Let $x_{i}=\displaystyle \max_{v_{k}\in V_{1}} x_{k}$ and $x_{j}=\displaystyle \max_{v_{k}\in V_{2}} x_{k}$. Then, considering the $i$th component of ${D(\mathcal{P}(\mathcal{D}_{2n}))}X=\rho_{1}X$, we get, 
\begin{eqnarray*}
\rho_{1} x_{i}=\displaystyle \sum_{v_{k}\in V_{1}}d_{i,k} x_{k} + \displaystyle \sum_{v_{k}\in V_{2}}d_{i,k} x_{k} \leq \displaystyle \sum_{v_{k}\in V_{1}}d_{i,k} x_{i} + \displaystyle \sum_{v_{k}\in V_{2}}d_{i,k} x_{j}\le T_{max}x_{i}+2nx_{j}
\end{eqnarray*}
\begin{equation}\label{inequ1dD}
    \implies (\rho_{1}-T_{max})x_{i}\leq 2nx_{j}
\end{equation}
Also, considering the $j$th component of ${D(\mathcal{P}(\mathcal{D}_{2n}))}X=\rho_{1} X$, we get, \begin{eqnarray*}
    \rho_{1} x_{j}=\displaystyle \sum_{v_{k}\in V_{1}}d_{i,k} x_{k} + \displaystyle \sum_{v_{k}\in V_{2}}d_{i,k} x_{k} \leq \displaystyle \sum_{v_{k}\in V_{1}}d_{i,k} x_{i} + \displaystyle \sum_{v_{k}\in V_{2}}d_{i,k} x_{j} \leq (2n-1)x_{i}+ 2(n-1)x_{j}
\end{eqnarray*}
\begin{equation}\label{inequ2dD}
    \implies (\rho_{1}-2n+2)x_{j}\leq (2n-1)x_{i}
\end{equation}
Since $x_{i},~ x_{j}$ are positive, from inequality (\ref{inequ1dD}) \& (\ref{inequ2dD}), we get $(\rho_{1}-T_{max})(\rho_{1}-2n+2) \leq 2n(2n-1).$
\begin{equation*}
    \implies \rho_{1}^{2}-(T_{max}+2n-2)\rho_{1}+\left[2(n-1)T_{max}-2n(2n-1)\right]\leq 0
\end{equation*}
As $\rho_{1} >0$
\begin{equation*}
\rho_{1}(\mathcal{P}(\mathcal{D}_{2n})) \leq \frac{1}{2}\left[(T_{max}+2n-2) + \sqrt{(T_{max}-2n+2)^{2}+8n(2n-1)}\right].
\end{equation*}
\end{proof}
Next, we provide the distance spectral radius bounds of the power graph of $\mathcal{Q}_{4n}$.
\begin{theorem}\label{dQ4n}
For any integer $n\ge 2$, the distance spectral radius $\rho_{1}(\mathcal{P}(\mathcal{Q}_{4n}))$ of $\mathcal{P}(\mathcal{Q}_{4n})$ satisfies the following: $$ \rho_{1}(\mathcal{P}(\mathcal{Q}_{4n})) \ge \frac{1}{2}\left[(T_{avg}+4n-3) + \sqrt{(T_{avg}-4n+3)^{2}+16(2n-1)^{2}}\right] $$ and $$ \rho_{1}(\mathcal{P}(\mathcal{Q}_{4n}))\le \frac{1}{2}\left[(T_{max}+4n-3) + \sqrt{(T_{max}-4n+3)^{2}+32n(2n-1)}\right]$$ where $T_{avg}$ and $T_{max}$ are the average and maximum transmissions of the graph $\mathcal{P}(C_{2n})$.
\end{theorem}
\begin{proof}  We are following the same vertex labeling and ordering as in Theorem \ref{aQ4n}. So the distance matrix of $\mathcal{P}(\mathcal{Q}_{4n})$ with rows and columns indexed by  $\{ V_{1}, V_{2}\}$ as the block matrix
    \begin{eqnarray*}\displaystyle {D(\mathcal{P}(\mathcal{Q}_{4n}))}&=&{\begin{bmatrix}  {D(\mathcal{P}(C_{2n}))} &{R_{2n}}\\
 {R^{T}_{2n}} &{S_{2n}} 
 \end{bmatrix}},
 \end{eqnarray*}
 where $R_{2n}$ is the block matrix $\begin{bmatrix}
     {J_{2 \times 2n}}\\
     {2J_{(2n-2) \times 2n}} 
 \end{bmatrix}$

 and $S_{2n}$ is the block matrix $\begin{bmatrix}
     {2(J_{n}-I_{n})} &{2J_{n}-I_{n}}\\
     {2J_{n}-I_{n}} &{2(J_{n}-I_{n})}
 \end{bmatrix}.$
 Then the quotient matrix of $D(\mathcal{P}(\mathcal{Q}_{4n}))$ with respect to the partition $\pi=\{ V_{1}, V_{2}\}$ is 
\begin{eqnarray*}\displaystyle {Q_{D}}&=&{\begin{bmatrix} 
 {T_{avg}} &{4n-2}\\
 {4n-2} &{4n-3} 
 \end{bmatrix}};
\end{eqnarray*} where $T_{avg}$ is the average transmission of $\mathcal{P}(C_{2n})$.
The spectral radius of $Q_{D}$ is \\ $\mu_{1}({Q_{D}})=\frac{1}{2}\left[(T_{avg} + 4n -3)+ \sqrt{(T_{avg}-4n+3)^{2}+16(2n-1)^{2}}\right]$.
So, by the Lemma \ref{quotient}, we get $\rho_{1}(\mathcal{P}(\mathcal{Q}_{4n}))\ge \mu_{1}({Q_{D}}).$
Therefore, $\rho_{1}(\mathcal{P}(\mathcal{Q}_{4n}))\ge \frac{1}{2}\left[(T_{avg} + 4n -3)+ \sqrt{(T_{avg}-4n+3)^{2}+16(2n-1)^{2}}\right].$

Again, $ D(\mathcal{P}(\mathcal{Q}_{4n}))$ is a non-negative and irreducible matrix. So, by the Perron-Frobenius Theorem, the spectral radius $\rho_{1}$ is simple and has a positive eigenvector, say $X=[x_{1}, x_{2}, \cdots, x_{4n}]^{T}.$\\
Let $x_{i}=\displaystyle \max_{v_{k}\in V_{1}} x_{k}$ and $x_{j}=\displaystyle \max_{v_{k}\in V_{2}} x_{k}$. Then, considering the $i$th component of ${D(\mathcal{P}(\mathcal{Q}_{4n}))}X=\rho_{1}X$, we get, 
\begin{eqnarray*}
\rho_{1} x_{i}=\displaystyle \sum_{v_{k}\in V_{1}}d_{i,k} x_{k} + \displaystyle \sum_{v_{k}\in V_{2}}d_{i,k} x_{k} \leq \displaystyle \sum_{v_{k}\in V_{1}}d_{i,k} x_{i} + \displaystyle \sum_{v_{k}\in V_{2}}d_{i,k} x_{j}\le T_{max}x_{i}+4nx_{j}
\end{eqnarray*}
\begin{equation}\label{inequ1dQ}
    \implies (\rho_{1}-T_{max})x_{i}\leq 4nx_{j}
\end{equation}
Also, considering the $j$th component of ${D(\mathcal{P}(\mathcal{Q}_{4n}))}X=\rho_{1} X$, we get, \begin{eqnarray*}
    \rho_{1} x_{j}=\displaystyle \sum_{v_{k}\in V_{1}}d_{i,k} x_{k} + \displaystyle \sum_{v_{k}\in V_{2}}d_{i,k} x_{k} \leq \displaystyle \sum_{v_{k}\in V_{1}}d_{i,k} x_{i} + \displaystyle \sum_{v_{k}\in V_{2}}d_{i,k} x_{j} \leq (4n-2)x_{i}+ (4n-3)x_{j}
\end{eqnarray*}
\begin{equation}\label{inequ2dQ}
    \implies (\rho_{1}-4n+3)x_{j}\leq (4n-2)x_{i}
\end{equation}
Since $x_{i},~ x_{j}$ are positive, from inequality (\ref{inequ1dQ}) \& (\ref{inequ2dQ}), we get $(\rho_{1}-T_{max})(\rho_{1}-4n+3) \leq 4n(4n-2).$
\begin{equation*}
    \implies \rho_{1}^{2}-(T_{max}+4n-3)\rho_{1}+[(4n-3)T_{max}-4n(4n-2)]\leq 0
\end{equation*}
As $\rho_{1} >0$
\begin{equation*}
\rho_{1}(\mathcal{P}(\mathcal{Q}_{4n})) \leq \frac{1}{2}\left[(T_{max}+4n-3) + \sqrt{(T_{max}-4n+3)^{2}+32n(2n-1)}\right].
\end{equation*}
Hence complete the proof.
\end{proof}

Next, we get the distance spectral radius of the family of the cyclic group of order $pq$, where $p$ and $q$ are primes, as in terms of the largest root of a cubic equation.
\begin{theorem}
    Let $\mathbb{G}=\mathbb{Z}_{pq}$, where $p$ and $q$ are primes, then the distance spectral radius of $D(\mathcal{P}(\mathbb{Z}_{pq}))$ is the largest root of the equation $x^{3}-(pq-3)x^{2} -(5pq-3p-3q)x+(p^{2}q^{2}-2p^{2}q-2pq^{2}+p^{2}+pq+q^{2})=0$.
\end{theorem}
\begin{proof}
  We consider the same vertex labeling and ordering as in Theorem \ref{Zpq}. As we observed that none of the vertices in $\mathcal{P}(V_{2})$ is adjacent to any vertex in $\mathcal{P}(V_{3})$ and graph diameter is $2$. So, the distance from any vertices in $\mathcal{P}(V_{2})$ to the every vertices in $\mathcal{P}(V_{3})$ is $2$.
  So, we can write the distance matrix of the power graph $\mathcal{P}(\mathbb{Z}_{pq})$ with rows and columns indexed by $\{ V_{1}, V_{2}, V_{3}\}$:
   \begin{eqnarray*} \displaystyle {D(\mathcal{P}(\mathbb{Z}_{pq}))}&=&{\begin{bmatrix} 
    {J_{\phi(pq)+1}-I_{\phi(pq)+1}} & {J_{(\phi(pq)+1)\times (q-1)}} & {J_{(\phi(pq)+1)\times (p-1)}}\\
 {J_{(q-1)\times (\phi(pq)+1)}} & {J_{(q-1)}-I_{(q-1)}} & {2J_{(q-1)\times (p-1)}}\\
 {J_{(p-1) \times (\phi(pq)+1)}} & {2J_{(p-1)\times (q-1)}} & {J_{(p-1)}-I_{(p-1)}}
 \end{bmatrix}}.
 \end{eqnarray*}
 
 Then the quotient matrix of $D(\mathcal{P}(\mathbb{Z}_{pq}))$ with respect to the partition $\pi=\{ V_{1}, V_{2}, V_{3}\}$ is 
 \begin{eqnarray*}\displaystyle {Q_{D}^{\mathcal{P}(\mathbb{Z}_{pq})}}&=&{\begin{bmatrix} 
 {(p-1)(q-1)} &{q-1} &{p-1}\\
 {(p-1)(q-1)+1} &{q-2} &{2p-2}\\
 {(p-1)(q-1)+1} &{2q-2} &{p-2}
 \end{bmatrix}}.
\end{eqnarray*} It can be observed that the partition $\pi=\{ V_{1}, V_{2}, V_{3}\}$ is an equitable partition. Therefore, the spectral radius of the matrix  $D(\mathcal{P}(\mathbb{Z}_{pq}))$ is nothing but the the spectral radius of the matrix $Q_{D}^{\mathcal{P}(\mathbb{Z}_{pq})}$. The Perron-Frobenius Theorem ensures that the spectral radius of the matrix  $D(\mathcal{P}(\mathbb{Z}_{pq}))$ is the largest root of the characteristic polynomial of the matrix $Q_{D}^{\mathcal{P}(\mathbb{Z}_{pq})}$ and the required characteristic equation is det$(xI-Q_{D}^{\mathcal{P}(\mathbb{Z}_{pq})})=0$ i.e., $x^{3}-(pq-3)x^{2} -(5pq-3p-3q)x+(p^{2}q^{2}-2p^{2}q-2pq^{2}+p^{2}+pq+q^{2})=0$.
\end{proof}
\section*{Acknowledgement}
We would like to express our sincere gratitude to Dr. Anirban Bose for his invaluable participation in discussions. The first author acknowledges financial support from the ANRF-CRG India and  SRC, IIT Hyderabad. The third author acknowledges the financial support from the NBHM, India.
\section*{Declarations}
\textbf{Conflict of interest:} The authors state that they possess no conflicts of interest.


\end{document}